\documentclass{gtmon_a}  
\pdfoutput=1

\usepackage{pinlabel}

\proceedingstitle{Heegaard splittings of 3-manifolds (Haifa 2005)}
\conferencestart{10 July 2005}
\conferenceend{19 July 2005}
\conferencename{Heegaard splittings of 3-manifolds}
\conferencelocation{Haifa}

\editor{Cameron Gordon}
\givenname{Cameron}
\surname{Gordon}

\editor{Yoav Moriah}
\givenname{Yoav}
\surname{Moriah}

\title{Degeneration of Heegaard genus, a survey}

\author{David Bachman}
\givenname{David}
\surname{Bachman}
\address{Mathematics Department\\Pitzer College\\\newline
Claremont, CA 91711\\USA}
\email{bachman@pitzer.edu}

\author{Ryan Derby-Talbot}
\givenname{Ryan}
\surname{Derby-Talbot}
\address{Department of Mathematics\\The American 
University in Cairo\\\newline
113 Kasr El Aini St\\PO Box 2511\\Cairo 11511\\Egypt}
\email{rdt@aucegypt.edu}


\arxivreference{math.GT/0606383}

\volumenumber{12}
\issuenumber{}
\publicationyear{2007}
\papernumber{01}
\startpage{1}
\endpage{15}
\doi{}
\MR{}
\Zbl{}
\received{15 June 2006}
\revised{31 January 2007}
\accepted{8 February 2007}
\published{3 December 2007}
\publishedonline{3 December 2007}
\proposed{}
\seconded{}
\corresponding{}
\version{}

\makeatletter
\def\cnewtheorem#1[#2]#3{\newtheorem{#1}{#3}[section]
\expandafter\let\csname c@#1\endcsname\c@thm}

\theoremstyle{plain}
\newtheorem{thm}{Theorem}[section]
\cnewtheorem{cor}[thm]{Corollary}
\cnewtheorem{lem}[thm]{Lemma}
\cnewtheorem{conj}[thm]{Conjecture}

\theoremstyle{definition}
\cnewtheorem{define}[thm]{Definition}
\cnewtheorem{question}[thm]{Question}
\cnewtheorem{claim}[thm]{Claim}

\cnewtheorem{proposition}[thm]{Proposition}
\cnewtheorem{remark}[thm]{Remark}
\cnewtheorem{remarks}[thm]{Remarks}
\cnewtheorem{example}[thm]{Example}

\makeatother

\newcommand{\hs}{$V \cup_H W$ }

\begin{document}

\begin{abstract}  
We survey known (and unknown) results about the behavior of Heegaard
genus of 3-manifolds constructed via various gluings. The
constructions we consider are (1) gluing together two 3-manifolds with
incompressible boundary, (2) gluing together the boundary components
of surface$\times I$, and (3) gluing a handlebody to the
boundary of a 3-manifold. We detail those cases in which it is known
when the Heegaard genus is less than what is expected after gluing.
\end{abstract}

\begin{asciiabstract}  
We survey known (and unknown) results about the behavior of Heegaard
genus of 3-manifolds constructed via various gluings. The
constructions we consider are (1) gluing together two 3-manifolds with
incompressible boundary, (2) gluing together the boundary components
of surface times I, and (3) gluing a handlebody to the
boundary of a 3-manifold. We detail those cases in which it is known
when the Heegaard genus is less than what is expected after gluing.
\end{asciiabstract}

\maketitle


\section{Introduction}

In this paper we survey known (and unknown) results about the Heegaard genus of 3-manifolds constructed via some gluing map. In particular we will be concerned here with compact orientable 3-manifolds $M$ constructed in one of three ways.  For the first construction, let $X$ and $Y$ be compact, orientable, irreducible 3-manifolds each with a single boundary component homeomorphic to a closed orientable surface $F$,  and suppose that $\partial X$ and $\partial Y$ are essential in $X$ and $Y$. 

\medskip

\noindent {\bf Construction 1}\qua Glue $\partial X$ to $\partial Y$
via a map $\varphi$. We write $M=X \cup _F\mskip -0.5mu Y$ or $M = X \cup_{\varphi}\mskip -0.5muY$. 

\medskip

For the second construction we begin with $F \times I$. 

\medskip

\noindent {\bf Construction 2}\qua Glue $F \times \{0\}$ to $F \times \{1\}$ via a map $\varphi$. We write $M=F \times_{\varphi} S^1$.

\medskip

Finally, for the last construction considered here let $\mathcal H (F)$ denote the handlebody whose boundary is homeomorphic to $F$. 

\medskip

\noindent {\bf Construction 3}\qua Glue $\partial \mathcal H(F)$ to  $\partial X$ via a map $\varphi$. We write $M=X \cup \mathcal H(F)$ or $M=X \cup_{\varphi} \mathcal H(F)$. 

\medskip

The {\em Heegaard genus} of $M$, denoted $g(M)$, is the minimal value $g$ such that $M$ admits a Heegaard splitting of genus $g$. The genus of the surface $F$ is denoted $g(F)$. 

The following definition can be made for any compact, orientable 3-manifold $M$. For simplicity, we assume here that $M$ is closed. 

\begin{define} For $F$ a separating surface in $M$, let $X$ and $Y$ denote the components of $M$ cut along $F$. Let $V_X \cup_{H_X} W_X$ and $V_Y \cup_{H_Y} W_Y$ denote Heegaard splittings of $X$ and $Y$, respectively, such that $F \subset \partial V_X, \partial W_Y$. Then there exists a product neighborhood $F \times I$ of $F$ such that $V_X$ equals $F \times [0, \frac{1}{2}]$ with 1-handles attached along $F \times \{0\}$ and $W_Y$ equals $F \times [\frac{1}{2}, 1]$ with 1-handles attached along $F \times \{1\}$. Form a homeomorphism of $M$ by deforming $F \times I$ to $F \times \{\frac{1}{2}\}$ so that the disks of attachment of the 1-handles in $V_X$ and in $W_Y$ end up disjoint on $F\times \{\frac{1}{2} \} = F$. This yields compression bodies $V = V_Y \cup \{\mbox{1-handles in $V_X$}\}$, and $W = W_X \cup \{\mbox{1-handles in $W_Y$}\}$, giving a Heegaard splitting $V \cup_H W$. Such a splitting is called an {\em amalgamation along $F$}.
\end{define}

\begin{figure}[h]
\centering
\labellist\small
\pinlabel $W_X$ at 133 620
\pinlabel $F$ [r] at 62 530
\pinlabel $V_X$ at 139 555
\pinlabel $V_Y$ at 172 431
\pinlabel $W_Y$ at 139 503
\pinlabel $W$ at 391 593
\pinlabel $V$ at 391 454
\pinlabel $M$ at 139 386
\endlabellist
\includegraphics[width=3.5in]{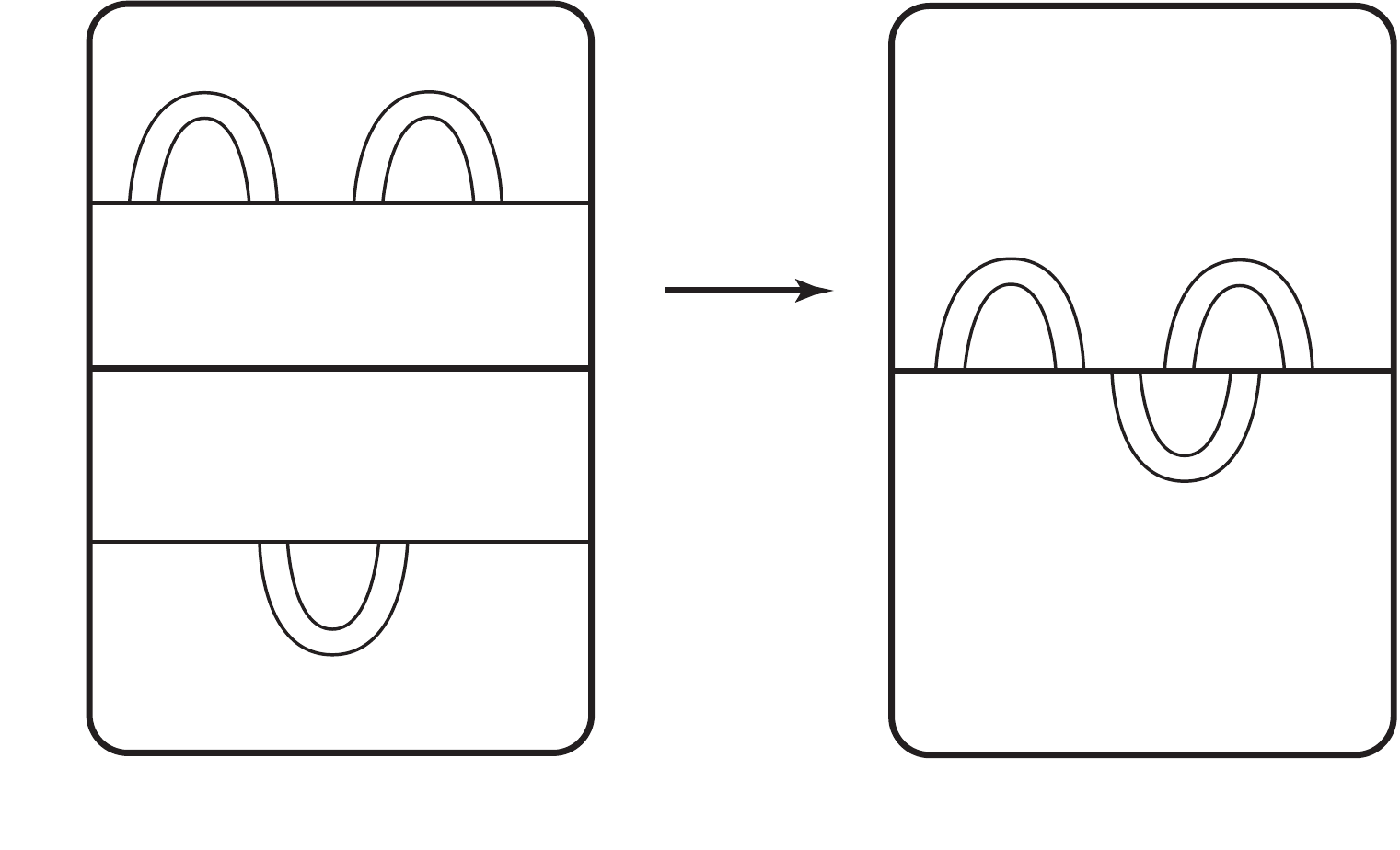}
\caption{A schematic for the construction of an amalgamation along $F$}
\end{figure}
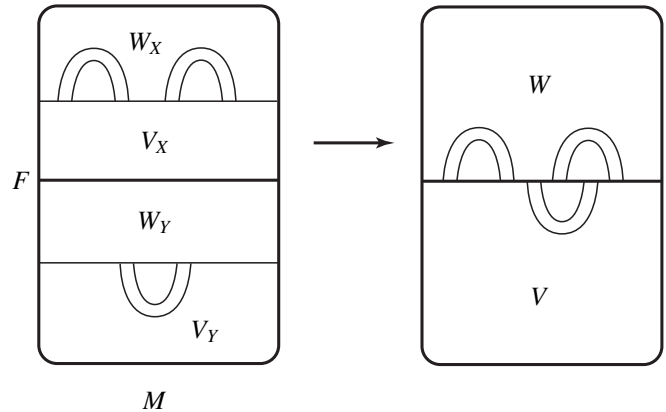

By amalgamating minimal genus Heegaard splittings along $F$ (or along two copies of $F$ if $M=F \times_{\varphi} S^1$) one obtains the following inequalities in each of the three constructions:
\medskip

\begin{enumerate}
	\item $g(M) \leq g(X) + g(Y) - g(F)$
	\item $g(M) \leq 2g(F) + 1$
	\item $g(M) \leq g(X)$ 
\end{enumerate}
\medskip

Our discussion will be concerned with when each of the inequalities is either strict (``degeneration is possible'') or is in fact an equality (``no degeneration''). 

If a Heegaard splitting of a 3-manifold is of minimal genus, then it is unstabilized. One way that degeneration of Heegaard genus can occur in the above situations is that the amalgamation of minimal genus splittings results in a stabilized splitting. This leads to the following natural question:

\begin{question}
\label{ques:stabilization}
When is the amalgamation of unstabilized Heegaard splittings unstabilized?
\end{question}

The paper is organized as follows. In Sections 2 -- 4 we discuss the
issue of degeneration of Heegaard genus under the three types of
gluing mentioned above, concluding each section with a discussion of
\fullref{ques:stabilization}. These sections are organized by the
genus of $F$ as the results tend to be based more on genus than on
which of the three constructions we are considering. The results on
Heegaard genus degeneration from these sections are summarized in a
table at the end of Section 4. In Section 5 we review known results
bounding how much degeneration of Heegaard genus can occur. For basic
definitions and notions related to Heegaard splittings, see
Scharlemann \cite{Scharlemann}.

An interesting property of Heegaard genus is that it provides an upper bound on the rank of the fundamental group of a 3-manifold $M$, since a genus $g$ Heegaard splitting of $M$ can be used to construct a presentation of $\pi_1(M)$ with $g$ generators. There are several results about the degeneration of rank of the fundamental group of manifolds formed via some gluing map, many analogous to those stated in this paper about Heegaard genus. While the analogy between Heegaard genus and rank is interesting and worth mentioning, we do not attempt to include a detailed discussion of it here. 

\section{Sphere gluings}

\subsection{No degeneration}

We begin by showing that Heegaard genus does not degenerate in any of the three constructions when $F$ is a sphere. The case that $M = X \cup_{S^2} Y$ is implied by the following classic result of Haken. 

\begin{thm}[Haken's Lemma \cite{Haken}]
Let \hs be a Heegaard splitting of $M$ and suppose that $M$ contains an essential 2-sphere. Then there is an essential 2-sphere $F$ such that $H \cap F$ is a single simple closed curve essential on $H$.  
\end{thm}

If the sphere $F$ is separating so that $M$ is the connected sum of two 3-manifolds $X$ and $Y$, then Haken's Lemma implies that a Heegaard splitting of $M$ is obtained from the ``connected sum'' of Heegaard splittings of $X$ and $Y$. In particular:  

\begin{cor}
\label{Haken}
If $M=X \# Y$, then $g(M) = g(X) + g(Y)$.
\end{cor}

\begin{proof}
Given Heegaard splittings of $X$ and $Y$ let $B_X$ and $B_Y$ be open embedded 3-balls in $X$ and $Y$, respectively, each intersecting the respective Heegaard surface in an open equatorial disk. Form the connected sum of $X$ and $Y$ by gluing $X - B_X$ to $Y - B_Y$ so that the component of the Heegaard surface in $X - B_X$ meets the component of the Heegaard surface in $Y - B_Y$ (note that there are two ways to do this, yielding possibly non-isotopic splittings). This yields a Heegaard splitting of $M$, implying that $g(M) \leq g(X) + g(Y)$.

For the reverse inequality, assume that $X$ and $Y$ are irreducible. Consider a Heegaard splitting \hs of $M$ of minimal genus. Since $M$ contains an essential sphere, by Haken's Lemma it contains one meeting $H$ in an essential simple closed curve. Since $X$ and $Y$ are irreducible, $M$ contains a unique essential sphere, hence $H$ intersects the connect sum sphere $S^2$ in a simple closed curve. Cutting along this sphere and filling the resulting boundary components with 3-balls containing equatorial disks, we obtain Heegaard splittings of $X$ and $Y$. If one of these splittings is not of minimal genus, then as above it could have been used to form a Heegaard splitting of $M$ of smaller genus, a contradiction. Thus the splittings of $X$ and $Y$ must be of minimal genus. This implies $g(M) = g(X) + g(Y)$.

To remove the assumption that $X$ and $Y$ are irreducible, use Milnor's result on the uniqueness of prime decomposition for 3-manifolds \cite{Milnor} and proceed by induction. 
\end{proof}
\medskip

If $M=S^2 \times_{\varphi} S^1$ is orientable then there is only one possibility for the map $\varphi$, which implies $M=S^2 \times S^1$. The Heegaard genus of this manifold is one.
\medskip

Finally, if $M$ has a sphere boundary and we glue a genus zero handle-body (a 3-ball) to it then the Heegaard genus does not change. Thus in all three constructions, Heegaard genus does not degenerate when $F$ is a sphere.

\subsection{Stabilization and connected sum}
In the case that $F$ is a sphere \fullref{ques:stabilization} was
originally asked by C\,McA Gordon, who conjectured that the connected
sum of unstabilized Heegaard splittings is never stabilized
\cite{kirby:97}.  A proof of this conjecture has been announced
independently by the first author and Qiu.

\begin{thm}{\rm\cite{Bachman,Qiu}}\qua
Let $H_X$ and $H_Y$ be unstabilized Heegaard surfaces in $X$ and $Y$, respectively. Then $H_X \# H_Y$ is an unstabilized Heegaard splitting surface in $X \# Y$. 
\end{thm}

The splitting surface $H_X \# H_Y$ is defined as in the proof of \fullref{Haken}. 

\section{Torus Gluing}

\subsection{Degeneration is possible}

Unlike the sphere case, there are examples where Heegaard genus degenerates when $F$ is a torus. The following result of Schultens and Weidmann shows that the amount of degeneration can be arbitrarily large.

\begin{thm}{\rm\cite{SchultensWeidmann}}\qua
\label{SchultensWeidmann}
Let $n$ be a positive integer. Then there exist manifolds $M_n = X_n \cup_{T^2} Y_n$ such that $$g(M_n) \leq g(X_n) + g(Y_n) - n.$$ 
\end{thm}

They in fact construct examples of unstabilized Heegaard splittings of $X_n$ and $Y_n$ such that the resulting amalgamated Heegaard splitting of $M_n$ can be destabilized $n$ times. 
\medskip

If $M = T^2 \times_{\varphi} S^1$ is a torus bundle, degeneration of Heegaard genus is also possible. Taking two genus 2 Heegaard splittings of $T^2 \times I$ and amalgamating gives a Heegaard splitting of genus 3 of $M$, implying that $g(M) \leq 3$. Cooper and Scharlemann have characterized precisely which solvmanifolds have $g(M)=2$.

\begin{thm}{\rm\cite{CooperScharlemann}}\qua
\label{CooperScharlemann}
Let $M = T^2 \times_{\varphi} S^1$ be a solvmanifold, and suppose the monodromy $\varphi$ can be expressed as $$\left( \begin{array}{cc} \pm m & -1 \\ 1 & 0\\ \end{array} \right).$$ If $m \geq 3$, then $g(M) = 2$. 
\end{thm}

Moreover, they show that there are precisely two genus 2 splittings if
$m = 3$, and only one genus 2 splitting if $m \geq 4$. It should be
noted that manifolds which are torus bundles but not solvmanifolds,
namely flat manifolds and nilmanifolds, have well understood Heegaard
splittings as they admit Seifert fibrations.  (See for example Moriah
and Schultens \cite{Moriah-Schultens} or Sedgwick \cite{Sedgwick} for
results on Heegaard genus of Seifert fibered spaces.)
\medskip

Finally, consider a 3-manifold obtained by gluing a solid torus to a
manifold $X$ with torus boundary, ie, via Dehn filling. In this case, as in the case of gluing two manifolds along a torus, Heegaard genus can degenerate by an arbitrary amount. 

\begin{example} 
\label{ex:tunnelnumber}
Let $X$ be the complement of a tunnel number $n$ knot. Perform trivial Dehn filling on $\partial X$ to obtain $S^3$. As the Heegaard genus of $S^3$ is $0$, it follows that for any $n$ there are manifolds $X$ such that
$$g(X \cup \mathcal H(F)) = g(X) - n.$$
\end{example}

Another way Heegaard genus can degenerate under Dehn filling is the following situation. 
\begin{example}
\label{ex:destabilizationlines}
Suppose that  $X$ has a single torus boundary component $T$ and let \hs be a Heegaard splitting of $X$. Assume that $T$ is contained in $V$, so $V$ is a compression body. Then there exists an essential disk $D'$ in $V$ such that $V$ cut along $D'$ contains a component $U$ homeomorphic to $T \times I$. Suppose there exists an essential disk $D$ in $W$ such that $\partial D$ meets the boundary of $U$ in a single arc $\delta'$. The endpoints of $\delta'$ lie on $\partial D'$. Let $\beta'$ be a loop on $\partial U$ composed of $\delta'$ and a properly embedded arc in $D'$. Then there is a loop $\beta$ on $T$ and an essential annulus in $U$ whose boundary components are $\beta$ and $\beta'$. Suppose $\alpha$ is a slope on $T$ meeting $\beta$ in a single point (there is an infinite number of such slopes). Attaching a solid torus $\mathcal H (F)$ by gluing a meridian disk to $\alpha$ makes $H$ a stabilized Heegaard surface in the resulting 3-manifold. These slopes correspond to a {\em destabilization line} in the Dehn filling space of $X$ (the {\em Dehn filling space of $X$} is the set of all 3-manifolds obtained by Dehn filling $X$). Thus in these situations, $$g(X \cup_{\alpha} \mathcal H(F)) \leq g(X) - 1.$$
\end{example}

\subsection{Sufficiently complicated torus gluings}
\label{sec:torusgluings}

Despite the fact that Heegaard genus can degenerate when gluing along
a torus, the following results show that degeneration is in fact a
special phenomenon. Recall that $X$ and $Y$ are irreducible and each
has a single incompressible torus boundary component ($X$ and $Y$ are
called {\em knot manifolds} in the terminology of Bachman, Schleimer,
and Sedgwick \cite{BSS}).

\begin{thm}{\rm\cite{BSS}}\label{BSS}\qua
Suppose that $\varphi  \colon  \partial X \to \partial Y$ is a sufficiently complicated homeomorphism.  Then the manifold $M(\varphi) = X \cup_\varphi Y$ has no strongly irreducible Heegaard splittings.
\end{thm}

The term {\em sufficiently complicated} is given in Definition 4.2 in
\cite{BSS} and is a technical statement about the distance $\varphi$
maps curves on the torus (for example a suitably large power of an Anosov map is sufficiently complicated).

It can be shown using the above theorem that every Heegaard splitting of $M(\varphi)$ is an amalgamation along $\partial X$, implying that $$g(M(\varphi)) = g(X) + g(Y) - 1.$$
If $M = T^2 \times_{\varphi} S^1$ is a solvmanifold, Scharlemann and Cooper's analysis applies here as well.

\begin{thm}{\rm\cite{CooperScharlemann}}\label{CooperScharlemann2}\qua
If $M= T^2 \times_{\varphi} S^1$ is a solvmanifold with monodromy $\varphi$ that cannot be expressed in the form given in \fullref{CooperScharlemann}, then the minimal genus Heegaard splitting of $M$ has genus equal to $3$ and is unique up to isotopy. 
\end{thm} 

Finally, suppose $M = X \cup \mathcal H(F)$ is obtained by Dehn
filling. Above we gave examples where $g(X)$ degenerates by an
arbitrarily large amount upon Dehn filling, and where $g(X)$ can
degenerate by at least one for all fillings along slopes corresponding
to a destabilization line in the Dehn filling space of $X$. Following
work of Rieck and Rieck--Sedgwick, we see that with mild assumptions
on $X$ the above situations are the only possible ways for $g(X)$ to
degenerate and are not generic occurrences.

\begin{thm}{\rm\cite{Rieck,RieckSedgwick}}\label{Rieck}\qua
Let $X$ be an acylindrical manifold with incompressible torus boundary $T$. Then 
\begin{enumerate}
\item there are only finitely many slopes on $T$ for which $$g(X \cup \mathcal H(F)) \leq g(X) - 2,$$ 
\item there are only finitely many destabilization lines in the Dehn filling space of $X$ such that $$g(X\cup \mathcal H(F)) \leq g(X) - 1.$$
\end{enumerate}
In particular, there are an infinite number of manifolds $X \cup \mathcal H(F)$ such that $$g(X \cup \mathcal H(F)) = g(X).$$
\end{thm}

Moriah and Rubinstein initially proved a similar theorem for
negatively curved manifolds in \cite{MR}.  Rieck and Sedgwick have
proven a more general version of the above theorem for small manifolds
in \cite{RieckSedgwick2}, restricting greatly the possibilities for
discrepancies between Heegaard splittings of $X \cup \mathcal H(F)$
and $X$ of any (not necessarily minimal) genus.

\subsection{Stabilization and amalgamation along a torus}

In considering \fullref{ques:stabilization}, the result of Schultens and Weidmann given in \fullref{SchultensWeidmann} shows that for any $n$ there exist examples of unstabilized Heegaard splittings that can be amalgamated to give a splitting that destabilizes $n$ times. It seems, however, that this situation is special. 

\begin{conj}
\label{c:GordonTorus}
Let $M = X \cup_{T^2} Y$ where $X$ and $Y$ each have a single incompressible torus boundary component. There is a complexity on maps $\varphi \colon \partial X \to \partial Y$ and an integer $n(X,Y)$ such that if the complexity of $\varphi$ is greater than $n$ then the amalgamation of any unstabilized splittings of $X$ and $Y$ is unstabilized.
\end{conj}

\section{Higher genus gluings}
Although the results are similar we consider the case when $g(F) \geq 2$ separately because the techniques and the implications of the theorems are different. For example, in the previous section the conclusion of \fullref{BSS} is that when $g(F) = 1$ and the gluing map is ``sufficiently complicated'', then $M$ contains {\em no} strongly irreducible Heegaard splittings. The results presented in \fullref{sec:highergenus} show that when $g(F) \geq 2$ and the gluing map is ``sufficiently complicated'', then $M$ contains {\em no minimal genus} strongly irreducible Heegaard splittings. 

\subsection{Degeneration is possible}

As with the case when $F$ is a torus, it is also possible for Heegaard
genus to degenerate when $F$ is a surface of genus at least $2$ as
shown by the following result of Kobayashi, Qiu, Rieck and Wang.

\begin{thm}{\rm\cite{KQRW}}\label{KQRW}\qua 
There exists a 3-manifold $M$ containing connected, separating
incompressible surfaces $F_n$ of arbitrarily large genus such that
amalgamating two minimal genus Heegaard splittings of $X_n$ and $Y_n$
along $F_n$ yields a $g(F_n) - 3$ times stabilized Heegaard splitting
of $M$.
\end{thm}

As a consequence it follows that 
$$g(M) \leq g(X_n) + g(Y_n) - 2g(F_n) + 3.$$ An interesting aspect of
these examples is that this degeneration occurs in the same 3-manifold
$M$, ie, $M$ does not depend on $n$.
\medskip

For manifolds of the form $F \times_{\varphi} S^1$, degeneration is also possible. 

\begin{example}
\label{ex:SFS}
Let $M$ be a Seifert fibered space with base $B$ a sphere and
containing three exceptional fibers of multiplicities $n$, $2n$, $2n$,
where $n$ is an integer greater than $2$. Assume that the Euler number
of $M$ is 0, so that there is some horizontal surface $F$ in $M$ (see
for example Hatcher \cite[Proposition 2.2]{Hatcher}). By Jaco
\cite[Theorem VI.34]{Jaco}, $M = F \times_{\varphi} S^1$. Moreover,
the surface $F$ branch covers $B$ (a sphere) and by an Euler
characteristic argument yields the following equation (see
\cite{Hatcher}):
$$\chi(F) = m\chi(B) - m\left(\frac{2n-1}{2n} + \frac{2n-1}{2n} + \frac{n-1}{n}\right)$$ $$= 2m - m \frac{6n - 4}{2n}$$ where $m$ is the degree of the cover. As the least common multiple of the multiplicities of the fibers divides $m$, it follows that $2n$ divides $m$. Moreover, the assumption that $n \geq 3$ implies that $2 - (6n-4)/2n$ is negative. Thus
$$\chi(F) = 2m - m\frac{6n - 4}{2n}$$
$$\leq 4n - (6n -4)$$
$$=4-2n.$$ Taking $\chi(F) = 2 - 2g(F)$ and solving, we obtain $$g(F)
\geq n-1.$$ Thus, a Heegaard splitting of $M$ which is an amalgamation
along $F$ has genus at least $2n -1$. It is well known, however, that
$g(M) = 2$ (see for example \cite{Boileau}). Therefore, given an
integer $n \geq 3$ there is a 3-manifold $M$ such that $g(M) \leq
2g(F) + 1 -(2n-3)$, implying Heegaard genus can degenerate by an
arbitrary amount.
\end{example}

Finally, for manifolds of the form $X \cup_F \mathcal H(F)$, as before degeneration can occur. This can be seen by modifying Examples \ref{ex:tunnelnumber} and \ref{ex:destabilizationlines}. 

\begin{example}
\label{ex:tunnelnumberhighergenus}
Let $\mathcal H(F)$ be a knotted handlebody in $S^3$ whose complement $X$ has incompressible boundary. Then there is a way of gluing $\mathcal H (F)$ to $X$ such that the resulting manifold is $S^3$. Thus if $g(F) = n$, $$g(X \cup \mathcal H (F)) \leq g(X) - n.$$
\end{example}

\begin{example}
\label{ex:destabilizationLinesHigherGenus}
Suppose that $X$ has a single boundary component $F$, where $g(F) \ge 2$.  Let \hs be a minimal genus Heegaard splitting of $X$. Assume that $F$ is contained in $V$, so $V$ is a compression body. For each  loop $\alpha$ on $F$ one can find an essential annulus in $V$ which meets $F$ in $\alpha$ and meets $H$ in a loop $\alpha _H$. Now let $D$ be a compressing disk for $H$ in $W$ and suppose there is a loop $\alpha$ on $F$ such that $\alpha _H$ meets $\partial D$ in a point. Attaching a handlebody $\mathcal H (F)$ in such a way so that $\alpha$ now bounds a disk makes $H$ a stabilized Heegaard surface in the resulting 3-manifold.
\end{example}

\subsection{Sufficiently complicated higher genus gluings}
\label{sec:highergenus}

A {\em simple} 3-manifold is a 3-manifold that is compact, orientable, irreducible, atoroidal, acylindrical and has incompressible boundary. The following result of Lackenby shows that when two simple 3-manifolds are glued along a surface $F$ with $g(F) \geq 2$ via a ``sufficiently complicated'' map, then as in the torus case there is no degeneration of Heegaard genus. 

\begin{thm}{\rm\cite{Lackenby1}}\label{Lackenby1}\qua
Let $X$ and $Y$ be simple 3-manifolds, and let $h \colon \partial X \to F$
and $h' \colon F \to \partial Y$ be homeomorphisms with some connected
surface $F$ of genus at least two. Let $\psi \colon F \to F$ be a
pseudo-Anosov homeomorphism. Then, provided $|n|$ is sufficiently
large,
\[g(X \cup_{h'\psi^n h} Y)=g(X)+g(Y)-g(F).\]
Furthermore, any minimal genus Heegaard splitting for $X
\cup_{h'\psi^n h} Y$ is obtained from splittings of $X$ and $Y$ by
amalgamation, and hence is weakly reducible.
\end{thm}

The intuition behind the proof of Lackenby's theorem is as
follows. When the map $\varphi$ is sufficiently complicated then
geometrically $M$ has a ``long neck" region homeomorphic to $F \times
(0,1)$. A result of Pitts and Rubinstein \cite{PittsRubinstein} implies that a strongly irreducible Heegaard splitting surface $H$ is isotopic to a minimal surface or to two copies of a double cover of a non-orientable minimal surface attached by a tube. In either case, if such a surface passes through the long neck region then it must have large area. By the Gauss-Bonnet theorem this implies $H$ has large genus. The conclusion is that if the map $\varphi$ is complicated enough then any strongly irreducible Heegaard splitting has genus higher than the genus of an amalgamated splitting. From here it is not difficult to show that any splitting (strongly irreducible or not) which is not an amalgamation of splittings of $X$ and $Y$ is not minimal genus. 

Souto has generalized this technique using the notion of distance in the curve complex. 

\begin{thm}{\rm\cite{Souto}}\label{Souto}\qua
Let $X$ and $Y$ be simple 3-manifolds and suppose $\partial X$ and
$\partial Y$ are connected and homeomorphic with genus at least
two. Fix an essential simple closed curve $\alpha \subset \partial X$
and $\alpha ' \subset \partial Y$. Then there is a constant $n_0$ such
that every minimal genus Heegaard splitting of $X \cup _ \psi Y$ is
constructed by amalgamating splittings of $X$ and $Y$ and hence
\[g(X \cup _\psi Y)=g (X ) + g (Y) - g (F)\] 
for every diffeomorphism $\psi: \partial X \to \partial Y$ with $d_{(\partial Y)} (\psi(\alpha), \alpha' ) \ge n_0$, where $d_{(\partial Y)} (\beta , \gamma)$ denotes the distance of essential simple closed curves $\beta$ and $\gamma$ in the curve complex 
of $\partial Y$. 
\end{thm}

Like Lackenby, Souto uses geometry to establish the above
result. T~Li has obtained a combinatorial proof of a similar theorem
\cite{Li}.
\medskip

Next suppose $M = F \times_{\varphi} S^1$. The following theorem of Lackenby indicates that, generically, the minimal genus Heegaard splittings of manifolds of the form $F \times_{\varphi} S^1$ are formed by amalgamating splittings of $F \times I$. 

\begin{thm}{\rm\cite{Lackenby2}}\label{Lackenby2}\qua
Let $M$ be a closed, orientable 3-manifold that fibers over the circle with pseudo-Anosov monodromy. Let $\{M_i \to M\}$ be the cyclic covers dual to the fiber. Then, for all but finitely many $i$, $M_i$ has an irreducible, weakly reducible, minimal genus Heegaard splitting.
\end{thm}

This implies that for all but finitely many $i$, $$g(M_i) = 2g(F)
+1.$$ Note that a stronger version of the above theorem has been
proved by Rubinstein \cite{Rubinstein}.  Also see a generalization by
Brittenham and Rieck \cite{RieckBrittenham}.

Bachman and Schleimer have improved this result using the notion of distance in the curve complex. Suppose that $M = F \times_{\varphi} S^1$ is formed using monodromy $\varphi \colon F \to F$. Define $d (\varphi)$ to be the minimum distance that $\varphi$ moves a vertex in the curve complex of $F$.

\begin{thm}{\rm\cite{BachmanSchleimer}}\label{BachmanSchleimer}\qua
Any Heegaard surface $H$ in $F \times_{\varphi} S^1$ with $-\chi(H) < d(\varphi)$ is an amalgamation of splittings of $F \times I$. 
\end{thm}

Finally, consider the case that $M$ is of the form $X \cup_F \mathcal
H(F)$. As a generalization of Thurston's hyperbolic Dehn surgery
theorem, Lackenby has shown in \cite{Lackenby3} that if $X$ is simple
and $\varphi \colon \partial X \to \partial \mathcal H(F)$ is
``sufficiently complicated'' then $X \cup_{\varphi} \mathcal H(F)$ is
irreducible, atoroidal, word hyperbolic and not Seifert
fibered. Lackenby then asks if the structure of the Heegaard
splittings of these manifolds can also be understood.

\begin{question}\cite{Lackenby3}\qua
How does Heegaard genus degenerate under handlebody gluing?
\end{question}

We pose the following conjecture as an answer to Lackenby's
question. This conjecture is similar in nature to \fullref{Rieck}.

\begin{conj}
Suppose that $X$ has a single boundary component $F$, where $g(F) \ge
2$.  Let \hs be a minimal genus Heegaard splitting of $X$. Assume that
$F$ is contained in $V$, so $V$ is a compression body. Let $\mathcal
W$ denote the set of vertices of the curve complex of $H$ that
correspond to the boundaries of disks in $W$. For each loop $\alpha$
on $F$ one can find an essential annulus in $V$ which meets $F$ in
$\alpha$ and meets $H$ in a loop $\alpha _H$. Now glue a handlebody
$\mathcal H(F)$ to $\partial X$. Let $\mathcal V_F$ denote the
vertices of the curve complex of $H$ defined as follows: if $\alpha$
bounds a disk in $\mathcal H(F)$ then $\alpha_H \in \mathcal V_F$. If
the distance between $\mathcal W$ and $\mathcal V_F$ is large enough
then $H$ is a minimal genus Heegaard splitting of $X \cup _F \mathcal
H(F)$.
\end{conj}

\subsection{Stabilization and amalgamation along a higher genus surface}

Again we consider \fullref{ques:stabilization}. As with the torus case, we have only a conjecture.

\begin{conj}
Let $M = X \cup_{\varphi} Y$ where $X$ and $Y$ each have a single incompressible boundary component of genus at least two. Then there is a complexity on maps $\varphi \colon \partial X \to \partial Y$ and an integer $n(X,Y,g)$ such that if the complexity of $\varphi$ is greater than $n$ the amalgamation of any unstabilized splittings of $X$ and $Y$ whose genera are less than $g$ is unstabilized.
\end{conj}

Note the subtle difference between this conjecture and \fullref{c:GordonTorus}. In \fullref{c:GordonTorus} we posit that if the gluing map is ``sufficiently complicated" then the amalgamation of {\it any} two unstabilized splittings is unstabilized. Here we conjecture that the same is true only if the splittings have low genus compared with the complexity of the gluing map.

\subsection{Table of degeneration of Heegaard genus}
In the following table we summarize the results of Sections 2 -- 4. The columns correspond to the genus of $F$ and the rows to the type of gluing used to construct $M$. We take ``D'' to mean ``degeneration is possible'', ``ND'' to mean ``no degeneration'', and ``NDSC'' to mean ``no degeneration if the gluing map $\varphi$ is sufficiently complicated'' in the appropriate contexts. In parentheses we provide the number of the theorem, corollary or example associated to the result.
\bigskip

\begin{center}
\begin{tabular}{|l||l|l|l|}
\hline
\ & $g(F)=0$ & $g(F) = 1$ & $g(F) \geq 2$\\ \hline \hline
$X \cup_F Y$ & ND \eqref{Haken}& D \eqref{SchultensWeidmann}& D \eqref{KQRW}\\ 
\ & \ & NDSC \eqref{BSS} & NDSC (\ref{Lackenby1}, \ref{Souto}) \\ \hline
$F \times_{\varphi} S^1$ & ND & D \eqref{CooperScharlemann} & D \eqref{ex:SFS} \\
\ & \ & NDSC \eqref{CooperScharlemann2}& NDSC (\ref{Lackenby2}, \ref{BachmanSchleimer}) \\ \hline
$X \cup_F \mathcal H(F)$ & ND & D (\ref{ex:tunnelnumber}, \ref{ex:destabilizationlines}) & D \eqref{ex:tunnelnumberhighergenus} \\
\ &  \ & NDSC \eqref{Rieck} &  NDSC ?? \\ \hline
\end{tabular}
\end{center}

\bigskip

\section{Lower bounds on the degeneration}

In the previous sections we discussed several situations in which
Heegaard genus can degenerate by an arbitrary amount. In this section
we state results that bound the amount by which Heegaard genus can
degenerate in terms of the genus of the gluing surface and the
Heegaard genera of the pieces. We will focus on the case that $M$ is
obtained by gluing $X$ and $Y$ together along a connected, orientable
surface $F$, ie, $M = X \cup_F Y$. Whereas some of the results on Heegaard genus degeneration in the previous sections are obtained by amalgamating unstabilized Heegaard splittings and getting stabilized splittings, the results in this section are obtained by finding lower bounds on the genus of the possible Heegaard splittings one can construct in a given manifold.  

As more restrictions are placed on the component manifolds $X$ and
$Y$, there are better known bounds. The least restrictive class of
manifolds was studied by Schultens \cite{Schultens}. Suppose $X$ and
$Y$ are irreducible 3-manifolds, and let $n_X$ and $n_Y$ denote the
number of non-parallel essential annuli that can be simultaneously
embedded in $X$ and $Y$, respectively. Then Schultens obtains the
bound
\[g(X \cup_F Y)\ge \frac{1}{5}(g(X)+g(Y)-8g(F)+11 - 4(n_X + n_Y)).\]
If, in addition, the manifolds $X$ and $Y$ are assumed the be
atoroidal and acylindrical, then previously Johannson \cite{Johannson}
had obtained the bound
\[g(X \cup_F Y)\ge \frac{1}{5}(g(X)+g(Y)-2g(F)).\]
Most recently the first author, in conjunction with Schleimer and
Sedgwick \cite{BSS}, added the restriction that the component
manifolds $X$ and $Y$ are {\it small} (ie, irreducible and
every incompressible surface is parallel to a boundary
component). This allowed them to obtain the bound
\[g(X \cup_F Y)\ge \frac{1}{2}(g(X)+g(Y)-2g(F)).\]

\bibliographystyle{gtart}
\bibliography{link}

\end{document}